\pdfoutput=1

\documentclass[twocolumn]{svjour3mod}

\input{arxiv_eigen_shared}

\DeclareMathOperator{\arctantwo}{arctan2}
\DeclareMathOperator*{\diagoper}{\text{diag}}

\DeclareMathOperator{\Corr}{\mathrm{Corr}}

\newcommand{\code}[1]{\texttt{#1}}

\newcommand{\ebox}[1]{\hbox{\ensuremath{#1}}}

\newcommand{\Identity}{I}
\newcommand{\inverse}{^{-1}}
\newcommand{\kronvector}{\boldsymbol{\delta}}

\newcommand{\ltwonorm}[1]{\Vert #1\Vert}

\newcommand{\orientedeigvec}{\mathcal{V}}

\newcommand{\realsspace}{\mathbb{R}}

\newcommand{\transpose}{^T}

\newcommand{\conv}{\ast}
\newcommand{\makevector}[1]{\mathbf{#1}}
\newcommand{\makevectorofsym}[1]{\boldsymbol{#1}}
\newcommand{\MarcenkoPastur}{Mar\v{c}enko\textendash\allowbreak{}Pastur\xspace}
\newcommand{\JamesStein}{James\textendash\allowbreak{}Stein\xspace}
\newcommand{\LedoitWolf}{Ledoit\textendash\allowbreak{}Wolf\xspace}
\newcommand{\makeestimate}[1]{\hat{#1}}
\newcommand{\etalwd}{et~al.\@\xspace}
\newcommand{\defeq}{\coloneqq}

\begin{document}

\title{A Consistently Oriented Basis for Eigenanalysis: Improved Directional Statistics}

\author{Jay Damask}
\authorrunning{Jay Damask}

\institute{
\at New York, NY. \\
\email{jaydamask@buell-lane-press.co} \\
orcid: 0000-0003-3288-2100
}

\date{Received: January 27, 2024}

\maketitle

\begin{abstract}
%
%

The algorithm derived in this article, which builds upon the original paper, takes a holistic view of the handedness of an orthonormal eigenvector matrix so as to transfer what would have been labeled as a reflection in the original algorithm into a rotation through a major arc in the new algorithm. In so doing, the angular wrap-around on the interval~$\pi$ that exists in the original is extended to a~$2\pi$ interval for primary rotations, which in turn provides clean directional statistics. The modified algorithm is detailed in this article and an empirical example is shown. The empirical example is analyzed in the context of random matrix theory, after which two methods are discussed to stabilize eigenvector pointing directions as they evolve in time. The \code{thucyd} Python package and source code, reported in the original paper, has been updated to include the new algorithm and is freely available.

  \keywords{Eigenvectors \and Eigenvalues \and Singular-value decomposition \and Directional statistics \and Random Matrices}
  \subclass{\\ 15A18 \and 65F15 \and 62H11 \and 15B52}
\end{abstract}


%
%

\section{Overview}
\label{eigen::sec: overview}

A model-free algorithm to consistently orient the eigenvectors returned from an \code{svd} or \code{eig} software call was detailed in an earlier article, and the Python package \code{thucyd} that was made freely available is now used in industry \cite{Damask:2020}. To consistently orient a matrix $V$ of eigenvectors is to replace $V$ with $\orientedeigvec$ such that the latter is only a pure rotation $R$ away from the identity matrix. To produce $\orientedeigvec$, the algorithm first sorts the eigenvectors by their associated eigenvalues in descending order, and then decides, one by one and from top to bottom, whether the sign of an eigenvector needs to be flipped. In the context of running linear regressions on the eigenbasis $\orientedeigvec$ over an evolving data stream, the original article demonstrated that the signs of regression weights $\makevectorofsym{\beta}$ are stabilized by regressing on $\orientedeigvec$ rather than $V$, which in turn leads to good interpretability of the evolving system. One of the consequences of the algorithm is that the rotation angles $\makevectorofsym{\theta}$ embedded within $V$ are calculated.%
\footnote{%
The angles are not unique, but given a fixed rotation order they can at least be standardized.
} %
As detailed then, the angles were limited to a \ebox{[-\pi/2, \pi/2]} interval and calculated using an \code{arcsin} approach. It was natural to extend the angular interval to \ebox{[-\pi, \pi]} by using an \code{arctan2} method, but it was shown that edge cases break the algorithm.

This article extends the original algorithm so that $N-1$ angles embedded in a space of $\realsspace^N$ are able to cover a full $2\pi$ interval. To do so requires the realization that only the first rotation upon entry into a new subspace of $V$ requires a $2\pi$ rotation range while the remaining angles in the subspace are inherently bound to a $\pi$ interval. This leads to a \emph{modified} \code{arctan2} algorithm that will not break on the aforementioned edge cases. As a consequence, the revised algorithm avoids angular wrap-around on the $\pi$ interval, a wrap-around that can be seen with the original. And yet, it also turns out that there is a tradeoff that is made between the original algorithm and its revision: sign flips in the original can be replaced by rotations on a $2\pi$ range (although not in a one-to-one manner), and as a consequence, the revised algorithm is better suited to problems that need full angular disambiguation while the original algorithm is better suited to problems like regression where a flip of a sign in $V$ has a knock-on effect to other parameters, such as $\makevectorofsym{\beta}$ above.

In either case, the orientation algorithm uses sequences of Givens rotations. The use of Givens rotations is far from new, and the original article cited Dash (2004) \cite{Dash:Polar:2004} who spoke of calculated the embedded angles in matrix $V$, and Golub and Van Loan (2013, 4th~ed.) \cite{Golub:2014}. Since that publication, I have become aware of the work of Cao \etalwd (2011) \cite{Cao:2011}, Tsay and Pourahmadi \cite{Tsay-Pourahmadi:2017}, and Potters and Bouchaud (2021) \cite{Potters-Bouchaud:2021}, all of whom use Givens rotations. However, the interest of Potters and Bouchaud is to calculate the Vandermonde determinant by perturbing the rotations, and the approaches of Cao and Tsay and company can lead to loss of positive definiteness. In the end, what is important here is the specification and solution to a system of transcendental equations, a sample of which appears in equation~(\ref{eigen::eq: trig vector equation def}). Golub and Van Loan come close to articulating this system with the second (unnumbered) equation in section 5.2.2, ``Householder QR,'' and the unnumbered equations in section 5.2.6, ``Hessenberg QR via Givens,'' (4$th$ edition), but they do not solve the equations. In effect, the algorithm in my original article solves the equations in one way and the revised algorithm detailed below solves the equations in another way. What is learned, in part, is that handedness, which is either left- or righthanded in 3D, does not generalize well in higher dimensions, but what can be said is whether the determinant of an orthonormal basis is $+1$, which indicates a pure rotation and an even number of reflections, or is $-1$, which indicates that there is an odd number of embedded reflections along with a rotation.

The body of this article reviews the original algorithm, details the revision, then applies the new tools to an empirical dataset, and connects the results to random matrix theory. Lastly, methods to stabilize eigenvectors as they evolve in time are demonstrated.

%
%

\section{Review of the Original Algorithm}
\label{eigen::sec: review of the original algorithm}

Following equation~(7) in the original paper, the central equation for the consistent orientation of an eigenvector matrix~$V$ is
\begin{equation}
  R\transpose V S = \Identity,
  \label{eigen::eq: RT V S = I representation def}
\end{equation}
where~$R$, $S$, and~$I$ are rotation, reflection, and identity matrices respectively. In the case where~$V$ has embedded reflections, $S$, which is defined by \ebox{S = \prod S_k} and
\begin{equation}
  S_k = \diagoper\left(1, 1, \ldots, s_k, \ldots, 1 \right),
  \label{eigen::eq: S_k def}
\end{equation}
where \ebox{s_k = \pm 1}, is designed to rectify them. With this design, the~$VS$ product produces
\begin{equation}
  \orientedeigvec = V S = R,
  \label{eigen::eq: V-orient = V S = R}
\end{equation}
where~$\orientedeigvec$ is the \emph{orientated} eigenvector matrix that is only a rotation~$R\transpose$ away from the identity. The algorithm that finds matrices \ebox{(R, S)} starts in the space of \ebox{V \in \realsspace^N} and visits its~$N-1$ subspaces so that all~$N$ eigenvectors are analyzed. (The spaces are indexed by the eigenvalues associated with the eigenvectors, sorted in descending order.) At each subspace level, beginning from the entire space, the algorithm determines whether an embedded reflection exists and if so, applies a reflection to the associated eigenvector. After this step, the rotation required to align that eigenvector to a constituent axis is calculated. The pattern that develops is
\begin{equation}
  V \rightarrow V S_1 \rightarrow R_1\transpose V S_1
                      \rightarrow R_1\transpose V S_1 S_2
                      \rightarrow \cdots ,
  \label{eigen::eq: program to attain R V S = I}
\end{equation}
which leads to the expansion of~(\ref{eigen::eq: RT V S = I representation def}) according to
\begin{equation}
  R_n\transpose \ldots R_2\transpose R_1\transpose\: V S_1 S_2 \ldots S_n = \Identity .
  \label{eigen::eq: representation expanded}
\end{equation}
For later convenience, let us define the partially oriented eigenvector $V_k$ matrix as
\begin{equation}
  V_k \defeq R_k\transpose \ldots R_1\transpose\: V S_1 \ldots S_k \;.
  \label{eigen::eq: Vk partial orientation}
\end{equation}
It is important to observe that reflections lie solely to the right of~$V$ whereas rotations lie solely to the left. Section~7.1 of the original paper considered the substitution of Householder reflections~$H$ for rotations~$R$, leading to \ebox{H\transpose V S}, and while this matrix product can be constructed to equal the identity matrix, it was shown that the reflections embedded in~$H$ (by its nature) scatter the carefully constructed reflections in~$S$, nullifying the interpretability of the solution.

\begin{figure*}[t!]
  \centering
  \includegraphics[width=\textwidth]{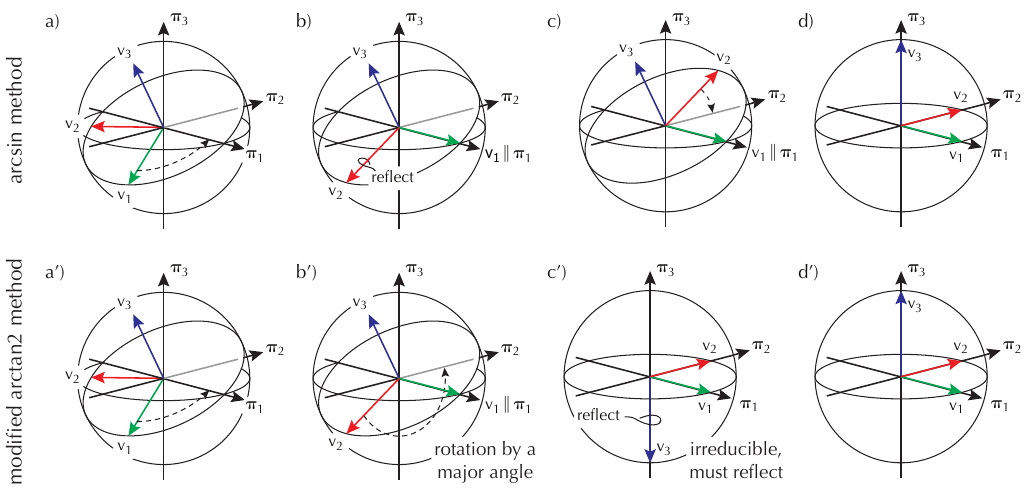}
  \caption{Two methods to orient a lefthanded basis in $\realsspace^3$ onto the righthanded identity matrix $\Identity$. The \emph{top} row shows the orientation sequence reported in the original article wherein $v_2$ as it appears in pane (b) is reflected in order to have it point in the direction of $\pi_2$. The \emph{bottom} rows show the new orientation sequence wherein rotation through a major angle replaces reflection in reducible subspaces while reflection is reserved for the last subspace, which is irreducible.
  }
  \label{eigen::fig: arctan2_reflect_vs_rotate}
\end{figure*}

To breathe life into the structure above, let us work through an example in~$\realsspace^3$ that is taken from figure~3 in the original paper and reproduced in the top row of figure~\ref{eigen::fig: arctan2_reflect_vs_rotate}. Pane~(a) shows a left-handed basis~$\makevector{v}$ that is to be oriented into the right-handed basis~$\makevectorofsym{\pi}$. At the first subspace (the entire space itself), the focus is on eigenvector~$v_1$ and constituent axis~$\pi_1$. As drawn, $v_1$ lies in the front hemisphere with respect to~$\pi_1$, which is to say, \ebox{v_1 \cdot \pi_1 > 0}, and so no reflection of~$v_1$ is required; thus \ebox{s_1 = 1}. Rotation~$R_1\transpose$ is then found to swing~$v_1$ onto~$\pi_1$; observe that~$R_1\transpose$ rotates the \emph{entire} basis~$\makevector{v}$, as pictured in pane~(b).

The next subspace, orthogonal to~$\pi_1$, is spanned by vectors \ebox{R_1\transpose[v_2, v_3]}. The algorithm begins again by constructing~$S_2$ first and then~$R_2\transpose$. As drawn, $R_1\transpose v_2$ lies in the back hemisphere with respect to~$\pi_2$, that is, \ebox{v_2 \cdot \pi_2 < 0}, which indicates the existence of an embedded reflection within $V$. The~$R_1\transpose$-rotated eigenvector~$v_2$ therefore requires reflection, so \ebox{s_2 = -1}, leading to the right-handed basis orientation seen in pane~(c). The orientation \ebox{R_1\transpose V S_1 S_2} requires a rotation $R_2\transpose$ about $\pi_1$ to align $v_2$ with $\pi_2$; this leads to pane~(d).

The first and second subspaces are \emph{reducible} because a rotation can align the associated eigenvector to a constituent axis, reducing the dimension of the nonaligned space by one. The last subspace, the third subspace in this example, is \emph{irreducible}: no axes remain about which a rotation can be applied; only reflection is possible. As drawn in pane~(c), $v_3$ lies in the front hemisphere with respect to $\pi_3$, so no reflection is required. To maintain symmetry in~(\ref{eigen::eq: representation expanded}), $R_3$ exists and is simply the identity matrix.

The pattern of working down the subspaces of $V$ can be seen in matrix form. Writing~$V$ of any dimension as
\begin{equation}
  V = \left(
    \begin{array}{ccc}
      & & \\
      & W_{(n)} & \\
      & &
    \end{array}
  \right),
  \label{eigen::eq: V = W_(n)}
\end{equation}
where~$W$ represents a working matrix and index~$n$ denotes the size of the subspace, the first two subspace reductions can be written as
\begin{equation}
  V_1 =
  \left(
    \begin{array}{cccc}
      1 & & & \\
        & \multicolumn{3}{c}{
            \multirow{3}{36pt}{ $W_{(n-1)}$ }} \\ \\ \\
    \end{array}
  \right)
  \;\;\text{and}\;\;
  V_2 =
  \left(
    \begin{array}{cccc}
      1 &   & & \\
        & 1 & & \\
        &   & \multicolumn{2}{c}{ W_{(n-2)} } \\
    \end{array}
  \right) .
  \label{eigen::eq: Rk W_(n-k) -> (1, 1, .., W_(n-k-1))}
\end{equation}
Observe that rotation by~$R_k\transpose$ leads to partial diagonalization of subspace $W_{(n-k-1)}$ such that a $1$ appears on the $k$th diagonal entry. This can be seen above as $V_1$ becomes \ebox{V_2 = R_2\transpose V_1 S_2}. The partial diagonalization continues for all reducible subspaces.

It is this pattern that leads to the governing equation to compute the angles of rotation. Let us label the entries of $W_{(n-k)}$ by the triplet \ebox{w_{n-k, i, j}}.\footnote{%
The original paper uses different indexing of the triplet. The indexing scheme here is clearer.
} %
The first two rotations require
\begin{equation*}
  R_1\transpose
    \left(
      \begin{array}{c}
        w_{n,1,1} \\
        w_{n,2,1} \\
        \vdots \\
        w_{n,n,1}
      \end{array}
    \right) =
    \left(
      \begin{array}{c}
        1 \\
        0 \\
        \vdots \\
        0
      \end{array}
    \right),
  \quad
  R_2\transpose
    \left(
      \begin{array}{c}
        0 \\
        w_{n-1,2,1} \\
        \vdots \\
        w_{n-1,2,n}
      \end{array}
    \right) =
    \left(
      \begin{array}{c}
        0 \\
        1 \\
        \vdots \\
        0
      \end{array}
    \right),
%
\end{equation*}
where the~$w$ columns naturally have unit $\ell^2$ norms. Rotations $R_k\transpose$ are then composed from a cascade of $n-k-1$ Givens rotations, each Givens rotation operating within a two-dimensional plane. For instance, \ebox{R_1 \in \realsspace^4} can be composed as
\begin{equation*}
  R_1(\theta_{1,2}, \theta_{1,3}, \theta_{1,4}) =
    R_{1,2}(\theta_{1,2})\: R_{1,3}(\theta_{1,3})\: R_{1,4}(\theta_{1,4}) .
\end{equation*}
The application of~$R_1\transpose$ to the first column of \ebox{W_{(4)}} leads to the composition shown figuratively as
\begin{equation}
  \underbrace{
      \left(
        \begin{array}{cccc}
          \centerdot & \centerdot &            & \\
          \centerdot & \centerdot &            & \\
                     &            & \centerdot & \\
                     &            &            & \centerdot
        \end{array}
      \right)
  }_{R_{1,2}}
  \;
  \underbrace{
      \left(
        \begin{array}{cccc}
          \centerdot &            & \centerdot & \\
                     & \centerdot &            & \\
          \centerdot &            & \centerdot & \\
                     &            &            & \centerdot
        \end{array}
      \right)
  }_{R_{1,3}}
  \;
  \underbrace{
      \left(
        \begin{array}{cccc}
          \centerdot &            &            & \centerdot \\
                     & \centerdot &            & \\
                     &            & \centerdot & \\
          \centerdot &            &            & \centerdot
        \end{array}
      \right)
  }_{R_{1,4}}
  \;
  \underbrace{
      \left(
        \begin{array}{c}
          \centerdot \\
                     \\
                     \\
                     \\
        \end{array}
      \right)
  }_{\kronvector_{1}}
  =
  \underbrace{
      \left(
        \begin{array}{c}
          \centerdot \\
          \centerdot \\
          \centerdot \\
          \centerdot
        \end{array}
      \right)
  }_{w_{4,,1}},
  \label{eigen::eq: Givens cascade in R4, dot notation}
\end{equation}
where, for instance,
\begin{equation*}
    R_{1, 3}\left(\theta_{1, 3}\right) =
    \left(
      \begin{array}{cccc}
        c_3  &     & -s_3  &  \\
             & 1   &       &  \\
        s_3  &     & c_3   &  \\
             &     &       & 1
      \end{array}
    \right)
    \longrightarrow
    \left(
      \begin{array}{cccc}
        \centerdot &            & \centerdot & \\
                   & \centerdot &            & \\
        \centerdot &            & \centerdot & \\
                   &            &            & \centerdot
      \end{array}
    \right),
  %
\end{equation*}
and \ebox{\{s_3|c_3\} = \{\sin|\cos\}(\theta_{1,3})}. Multiplying through the lefthand side of~(\ref{eigen::eq: Givens cascade in R4, dot notation}) and equating it to the righthand side leads to the system of transcendental equations
\begin{equation}
  \left(
    \begin{array}{c}
      c_2\: c_3 \: c_4 \\
      s_2\: c_3 \: c_4 \\
      s_3\: c_4        \\
      s_4
    \end{array}
  \right)
  =
  \left(
    \begin{array}{c}
      a_1 \\
      a_2 \\
      a_3 \\
      a_4
    \end{array}
  \right),
  \label{eigen::eq: trig vector equation def}
\end{equation}
where~$a_i$ is a simpler symbolic substitute for~$w_{4,,1}$. An alternative solution to this equation than the solutions found in the original article is central to this paper.

The original solution uses the \code{arcsin} method. Here, the entries in the columns are taken pairwise from the bottom up, leading to the system (again in $\realsspace^4$)
\begin{equation}
  \begin{aligned}
    \theta_{1,4} & = \arcsin \left( a_4 \right)   \\
    \theta_{1,3} & = \arcsin \left( a_3\, / \cos\theta_{1,4} \right)  \\
    \theta_{1,2} & = \arcsin \left( a_2\, / \left( \cos\theta_{1,4} \cos\theta_{1,3} \right)\right) .
  \end{aligned}
  \label{eigen::eq: arcsine solution to givens angles}
\end{equation}
Systems like this are composed for all reducible subspaces and the collection of angles $\makevectorofsym{\theta}$ are recorded. In all, there are \ebox{N(N-1)/2} angles embedded in \ebox{V\in\realsspace^N}.

Now, the issue with the \code{arcsin} method is that angles are minor, that is, those
restricted to a $\pi$ interval such as \ebox{\theta \in [-\pi/2, \pi/2]}. This in turn means that apparent the pointing direction of all~$N$ eigenvectors is limited to a hyper\-/hemisphere rather than a full hypersphere. In systems with a good mixture of states, the vectors tend to point toward the centers of their respective orthants, but in other systems, vectors can point near the equator of their respective hyper\-/hemispheres. In the presence of noise in an evolving system, vector wobble can lead to an angular wrap-around such that $v(t_2)$ and $v(t_1)$ appear on opposite sides of the hyper\-/hemisphere even when these vectors are in proximity on a hypersphere.

With this in mind, the original article proposed an \code{arctan2} method to span the major angular range of $2\pi$. However, it was shown that edge cases exist for this method that break the orientation algorithm, and so the \code{arctan2} method could not be used. The first version of the \code{thucyd} Python package was implemented with the \code{arcsin} method.

%
%

\section{Handedness in High Dimensions}
\label{eigen::sec: handedness in high dimensions}

This section revisits the analysis of equation~(\ref{eigen::eq: RT V S = I representation def}) and does so by considering the handedness of an orthogonal basis holistically. This approach contrasts with the ``opportunistic'' approach taken in the original article. The resulting \emph{modified} \code{arctan2} method is not simply a tweak to the existing algorithm but instead relies on insight about handedness and partially oriented eigenvectors.

\begin{figure*}[t!]
  \centering
  \includegraphics[width=125mm]{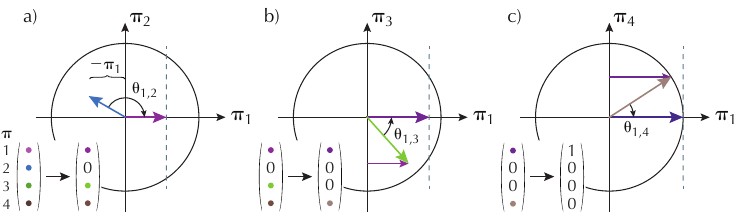}
  \caption{The orientation of a 4D vector via the three Givens rotations in~(\ref{eigen::eq: Givens cascade to orient a-vec to 1-vec}). The blue vector in pane (a) shows the orientation of \ebox{(a_1, a_2)} in the \ebox{(\pi_1, \pi_2)} plane. This vector is rotated into $\pi_1$, annihilating the $\pi_2$ component. Subsequent rotations, panes (b) and (c), annihilate the $\pi_3$ and $\pi_4$ components. Only in pane (a) does the vector point into the back hemisphere with respect to $\pi_1$.
  }
  \label{eigen::fig: arctan2_r4_rotations}
\end{figure*}

To begin, eigenvector matrix~$V$ and rotation matrix~$R$ are unitary in that
\begin{equation}
  V\transpose V = \Identity \quad\text{and}\quad R\transpose R = \Identity,
  \label{eigen::eq: VTV = I and RTR = I}
\end{equation}
and yet the individual matrices have different determinant possibilities:\,%
\footnote{%
The statement that \ebox{\det(V) = +1} on page four of the original article erroneously omits the potential for a negative sign.
} %
\begin{equation}
  \det(V) = \pm 1 \quad\text{and}\quad \det(R) = 1.
  \label{eigen::eq: }
\end{equation}
The reason for this difference is that~$V$ can contain embedded reflections along with its rotations, whereas $R$ is by nature a pure rotation. Therefore, $V$ belongs to the orthogonal group $O(N)$ while $R$ belongs to the special orthogonal group \ebox{SO(N) \subset O(N)}. The problem statement is then one of squeezing $V$ into $SO(N)$ so that the squeezed version of $V$ is only a pure rotation away from $\Identity$. To accomplish this, the auxiliary matrix~$S$ is introduced to impart pure reflections. Moving $R$ to the righthand side of~(\ref{eigen::eq: RT V S = I representation def}) and taking determinants leads to
\begin{equation}
  \det(V) \det(S) = \det(R) = 1.
  \label{eigen::eq: det(V) det(S) = det(R)}
\end{equation}
In the original scheme, the reflection matrix goes as
\begin{equation}
  S = \diagoper\left(\pm 1, \pm 1, \ldots, \pm 1\right),
  \label{eigen::eq: S = diag(pm 1, pm1, ...)}
\end{equation}
where the signs of the elements are determined as the subspaces are visited in descending order. The solution leads to an even number of negative entries (reflections) when \ebox{\det(V) = 1} and to an odd number of reflections otherwise, and no global constraint is placed on the number of elements of either sign. But in light of (\ref{eigen::eq: det(V) det(S) = det(R)}), it is clear that no reflections are inherently required when \ebox{\det(V) = 1}, and only a single reflection is necessary otherwise. To achieve the minimum reflection count, we need a workable method for major-angle rotation; this is the starting point for the revised algorithm.

Consider again the lefthanded $\realsspace^3$ basis shown in figure~\ref{eigen::fig: arctan2_reflect_vs_rotate}~\emph{top}. The diagonal entries were found to be \ebox{\diagoper(S) = \left(1, -1, 1\right)}, and indeed, being lefthanded, $V$ requires at least one reflection to become righthanded. But \emph{where} the reflection is introduced can be adjusted. The original algorithm is ``opportunistic'' because the hyper\-/hemisphere in which the $k$th eigenvector lies (after prior rotations) is forced via reflection (if required) to lie in the direction of constituent axis~$\pi_k$; the reflection is imposed by setting \ebox{s_k = -1}. The key way forward is to recognize that this reflection can be replaced by a rotation of $\pi$ when in a reducible subspace.

To illustrate this alternative approach, figure~\ref{eigen::fig: arctan2_reflect_vs_rotate}~\emph{bottom} again aligns a lefthanded $\realsspace^3$ basis. Panes (a) and (a$'$) are the same, since $v_1$ already points in the direction of~$\pi_1$, so no reflection is necessary; only a minor-angle rotation is required align the vectors. But while pane (b) reflects $v_2$ to point in the direction of $\pi_2$, the operation in pane (b$'$) \emph{rotates} $v_2$ through a major angle to align directly with $\pi_2$. Pane (c$'$) shows that $v_3$ is counterdirectional to $\pi_3$, and that there are no more rotations available with which to align these vectors\textemdash{}the system is irreducible. The only recourse is to reflect $v_3$ into $\pi_3$, and this is done to go from pane (c$'$) to (d$'$). The resulting reflection matrix is \ebox{\diagoper(S') = \left(1, 1, -1\right)}. Observe that the embedded reflection is now handled in the one irreducible dimension. There is in fact a pattern in how $S$ transforms into $S'$ as the subspaces are visited, and this is discussed below.

Moving on to higher dimensions, one way to draw an $N$-dimensional system is to plot the entries of a vector pairwise in a 2D plane across a series of panes. An example appears in figure~\ref{eigen::fig: arctan2_r4_rotations}. What is shown here is the evolution of a 4D vector as it is rotated into alignment with $\pi_1$ according to the rotations that appear in a reorganization of~(\ref{eigen::eq: Givens cascade in R4, dot notation}):
\begin{equation}
  \underbrace{
    \left(
      \begin{array}{cccc}
        \centerdot &            &            & \centerdot \\
                    & \centerdot &            & \\
                    &            & \centerdot & \\
        \centerdot &            &            & \centerdot
      \end{array}
    \right)
  }_{R_{1,4}\transpose}
  \underbrace{
    \left(
      \begin{array}{cccc}
        \centerdot &            & \centerdot & \\
                    & \centerdot &            & \\
        \centerdot &            & \centerdot & \\
                    &            &            & \centerdot
      \end{array}
    \right)
  }_{R_{1,3}\transpose}
  \underbrace{
    \left(
      \begin{array}{cccc}
        \centerdot & \centerdot &            & \\
        \centerdot & \centerdot &            & \\
                    &            & \centerdot & \\
                    &            &            & \centerdot
      \end{array}
    \right)
  }_{R_{1,2}\transpose}
  \left(
    \begin{array}{c}
      a_1 \\
      a_2 \\
      a_3 \\
      a_4
    \end{array}
  \right)
=
  \left(
    \begin{array}{c}
      1 \\
      0 \\
      0 \\
      0 \\
    \end{array}
  \right).
  \label{eigen::eq: Givens cascade to orient a-vec to 1-vec}
\end{equation}
Let us call the vector on the left $\makevector{a}$. As pictured in pane (a), $\makevector{a}$ points into the back hyper\-/hemisphere with respect to $\pi_1$. The original algorithm would reflect $\makevector{a}$ into the front hyper\-/hemisphere and then rotate, but here, rotation $R_{1,2}\transpose$ through major angle $\theta_{1,2}$ is used instead. The rotation is only within the \ebox{(\pi_1, \pi_2)} plane and applies only to the \ebox{(a_1, a_2)} entries. The figure inset shows that the $\pi_2$ component of \ebox{R_{1,2}\transpose\, \makevector{a}} vanishes. Next, the $R_{1,2}\transpose$-rotated \ebox{(a_1, a_3)} vector, shown in green in pane (b), is rotated by $R_{1,3}\transpose$ in the \ebox{(\pi_1, \pi_3)} plane so that the $\pi_3$ component is annihilated. Likewise, the $R_{1,3}\transpose R_{1,2}\transpose$-rotated \ebox{(a_1, a_4)} vector, the copper vector in pane (c), is rotated by $R_{1,4}\transpose$ so as to annihilate the $\pi_4$ component. To recapitulate the main point: a rotation through a major angle has substituted for a reflection and subsequent rotation through a minor angle.

Now, let us examine the signs of the entries. At the start, elements \ebox{(a_1, a_2)} can possess either sign: their vector direction can point into any quadrant in the \ebox{(\pi_1, \pi_2)} plane, and the \code{arctan2} function is powerful enough to determine the major angle $\theta_{1,2}$ required to align the vector to $\pi_1$. However, once the $a_1$ entry is rotated to align with $\pi_1$, subsequent values of the top element of the vector are guaranteed to be positive. As pane (b) shows, the third vector entry can be either positive or negative, but the angle $\theta_{1,3}$ is only ever minor. The same holds for $\theta_{1,4}$ in pane (c). Put another way, the vectors that lie in the \ebox{(\pi_1, \pi_3)} and \ebox{(\pi_1, \pi_4)} planes, prior to their orienting rotations, always point into the front hyper\-/hemisphere with respect to $\pi_1$. This observation lies at the heart of the modified \code{arctan2} method.

The modified \code{arctan2} solution to~(\ref{eigen::eq: trig vector equation def}), for four dimensions, is
\begin{equation}
  \begin{aligned}
    \theta_{1,2} & = \arctantwo\left( a_2, a_1 \right) \\
    \theta_{1,3} & = \arctantwo\left( a_3, \left| a_2 \csc \theta_2 \right| \right) \\
    \theta_{1,4} & = \arctantwo\left( a_4, \left| a_3 \csc \theta_3 \right| \right) ,
  \end{aligned}
  \label{eigen::eq: modified arctan2 solutions for givens angles}
\end{equation}
where the arguments follow a \code{arctan2(y, x)} call style. To avoid overflow, the latter equations can be rewritten in the form
\begin{equation}
  \theta_{\cdot,k} = \arctantwo\left( a_k \left| \sin \theta_{k-1} \right|, \left| a_{k-1} \right| \right) .
  \label{eigen::eq: modified arctan2 rearranged}
\end{equation}
In contrast to equation~(25) of the original article, only the first angle is calculated from all four quadrants; subsequent angles are calculated by enforcing the nonnegative nature of the vector as projected onto the $\pi_1$ axis. As a side effect, the catastrophic edge case discussed in the original article is avoided because a signed zero cannot appear in the \code{x} argument.

In the special case when the column vector~$\makevector{a}$ is sparse, the indexing in~(\ref{eigen::eq: modified arctan2 solutions for givens angles}) requires adjustment. To see this, consider \ebox{\makevector{a} = \left( a_1, a_2, 0, a_4 \right)\transpose}. Using~(\ref{eigen::eq: modified arctan2 solutions for givens angles}) as written, \ebox{a_3 = 0} leads to \ebox{\theta_{1,3} = 0}, which in turn leads to \ebox{\theta_{1,4} = \arctantwo(a_4 \cdot |0|, |0|) = 0} (from~(\ref{eigen::eq: modified arctan2 rearranged})). What has happened is that the algorithm, as written, has stalled, and the resulting ``solution'' to~(\ref{eigen::eq: trig vector equation def}) is incorrect. To correct for sparsity, let us first refer back to the original setup~(\ref{eigen::eq: Givens cascade to orient a-vec to 1-vec}): having \ebox{a_3 = 0} is tantamount to skipping the \ebox{R_{1,3}\transpose} rotation. With this in mind, all that needs to change is the indexing of the prior entry; that is,
\begin{equation}
  \theta_{\cdot,k} = \arctantwo\left( a_k \left| \sin \theta_{k-j} \right|, \left| a_{k-j} \right| \right) , \quad j \geq 1,
  \label{eigen::eq: modified arctan2 sparse}
\end{equation}
where \ebox{a_{k-j}} points to the first nonzero entry behind \ebox{a_k}. There is one more edge case when \ebox{a_2 = 0}. Here,
\begin{equation}
  \theta_{\cdot,k} = \arctantwo\left( a_k, \left| a_{1} \right| \right),
  \label{eigen::eq: modified arctan2 sparse a_2 = 0}
\end{equation}
where $k$ is the first nonzero entry in $\makevector{a}$ after $a_1$. The sine term is absent here, which reflects both its absence in the first row of~(\ref{eigen::eq: trig vector equation def}) and \ebox{\theta_{\cdot, j} = 0} for \ebox{j < k}.

All of the Givens rotation angles associated with the reducible subspaces of $V$ can be calculated in this way, leaving only a final potential reflection for the one irreducible subspace that remains. The central equation to orient $V$ into $\Identity$, equation~(\ref{eigen::eq: RT V S = I representation def}), is still
\begin{equation*}
  R_{\text{tan}}\transpose V S_{\text{tan}} = \Identity
\end{equation*}
in structure, but the angles $\makevectorofsym{\theta_\text{tan}}$ differ from the \code{arcsin} method, and the reflection matrix is simply
\begin{equation}
  S_{\text{tan}} = \diagoper\left( 1, 1, \ldots, \pm 1 \right).
  \label{eigen::eq: S for arctan2 method -- strict}
\end{equation}
In particular, $\makevectorofsym{\theta_\text{tan}}$ can now contain major angles for the first rotation of each new subspace.

To emphasize the point made in the original article, the sequence of Givens rotations to align a vector to a particular constituent axis is arbitrary and does not affect the ability to orient. Therefore, a convenient sequence and subsequent consistency are all that is required. It turns out that the Givens sequence for the \code{arcsin} method works just as well for the modified \code{arctan2} method.

\begin{figure}[t]
  \centering
  \includegraphics[width=84mm]{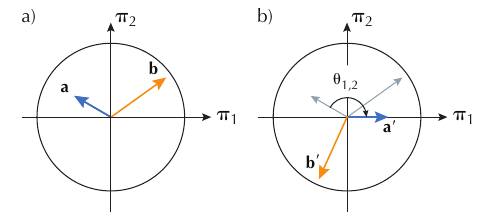}
  \caption{Rotation of a vector $\makevector{a}$ through a major angle, as projected onto the \ebox{(\pi_1, \pi_2)} plane, changes the hyper\-/hemisphere in which the vector points. The same rotation also reorients the \emph{next} vector in the system, here labeled $\makevector{b}$. This effect can be summarized as modifying two adjacent entries in reflection matrix $S$.
  }
  \label{eigen::fig: arctan2_coupled_rotation}
\end{figure}

The final topic of this section is the relationship between reflection matrix $S$ from the \code{arcsin} method, equation~(\ref{eigen::eq: S = diag(pm 1, pm1, ...)}), and the modified \code{arctan2} method, equation~(\ref{eigen::eq: S for arctan2 method -- strict}). One way to do this is to augment the $\makevector{a}$ vector in~(\ref{eigen::eq: Givens cascade to orient a-vec to 1-vec}) to write
\begin{equation}
  \underbrace{
    \left(
      \begin{array}{cccc}
        \centerdot &            &            & \centerdot \\
                    & \centerdot &            & \\
                    &            & \centerdot & \\
        \centerdot &            &            & \centerdot
      \end{array}
    \right)
  }_{R_{1,4}\transpose}
  \underbrace{
    \left(
      \begin{array}{cccc}
        \centerdot &            & \centerdot & \\
                    & \centerdot &            & \\
        \centerdot &            & \centerdot & \\
                    &            &            & \centerdot
      \end{array}
    \right)
  }_{R_{1,3}\transpose}
  \underbrace{
    \left(
      \begin{array}{cccc}
        \centerdot & \centerdot &            & \\
        \centerdot & \centerdot &            & \\
                    &            & \centerdot & \\
                    &            &            & \centerdot
      \end{array}
    \right)
  }_{R_{1,2}\transpose}
  \left(
    \begin{array}{cc}
      a_1 &  b_1 \\
      a_2 &  b_2 \\
      a_3 &  b_3 \\
      a_4 &  b_4
    \end{array}
  \right) .
  \label{eigen::eq: Givens cascade applied to vec-a and vec-b}
\end{equation}
Although \ebox{\makevector{a} \cdot \makevector{b} = 0} by construction, the entries that lie in the \ebox{(\pi_1, \pi_2)} plane are in general not orthogonal. An example is pictured in figure~\ref{eigen::fig: arctan2_coupled_rotation}~(a). Consider now the rotation by $R_{1,2}\transpose$ through a major angle. Although the two vectors, projected onto the \ebox{(\pi_1, \pi_2)} plane are not orthogonal, their angle is preserved under this rotation, as shown in pane~(b). And while $\makevector{a}$ and $\makevector{a}'$ point in opposite hyper\-/hemispheres, the \emph{same holds true} for vectors $\makevector{b}$ and $\makevector{b}'$ with respect to the hyper\-/hemisphere associated with axis $\pi_2$. This observation does not extend to further augmented vectors, say $\makevector{c}$ and $\makevector{d}$, because Givens rotations $R_{1,3}\transpose$ and $R_{1,4}\transpose$ subtend only minor angles. A rule can be deduced from this observation: a rotation through a major angle in subspace $k$ flips the sign of both $s_k$ and $s_{k+1}$. Reflection vectors $\diagoper(S_{\text{sin}})$ can now be ``untwisted,'' at least as a thought experiment.

To simplify the notation, let \ebox{\mathcal{S} = \diagoper(S)}. Then, from the example of figure~\ref{eigen::fig: arctan2_reflect_vs_rotate}, we found that $\mathcal{S_\text{sin}} = (1, -1, 1)$. Rather than applying a reflection in the second subspace, we can instead apply a major-angle rotation. This flips the sign of the second and third entries, yielding $\mathcal{S_\text{tan}} = (1, 1, -1)$. As another example, start with the system where $\mathcal{S_\text{sin}} = (-1, -1, -1, 1)$. Rotation in the first subspace by a major angle leads to $\mathcal{S} = (1, 1, -1, 1)$, after which only a minor rotation is needed in the second subspace. Another major-angle rotation in the third subspace gives $\mathcal{S_\text{tan}} = (1, 1, 1, -1)$, leaving a final reflection in the irreducible subspace. Observe that when \ebox{\det(V) = +1} and $\mathcal{S}_{\text{sin}}$ has an even number of negative signs, no reflection is needed for orientation using the modified \code{arctan2} method.

%
%

\section{Reference Implementation}
\label{eigen::sec: reference implementation}

The Python package \code{thucyd} is freely available from \textsc{PyPi} and \textsc{Conda-Forge} and can be used directly once installed. The source code is available on \textsc{GitLab}. See section~\ref{eigen::sec: declarations} for link details.

\SetAlCapSty{}
\begin{algorithm}[t!]
  {
    \DontPrintSemicolon
    \SetKw{KwStep}{step}
    \SetKwData{V}{V}
    \SetKwData{Vcol}{Vcol}
    \SetKwData{Vor}{Vor}
    \SetKwData{Vsort}{Vsort}
    \SetKwData{Vwork}{Vwork}
    \SetKwData{E}{E}
    \SetKwData{Eor}{Eor}
    \SetKwData{Esort}{Esort}
    \SetKwData{OrientToFirstOrthant}{OrientToFirstOrthant}
    \SetKwData{R}{R}
    \SetKwData{cursor}{cursor}
    \SetKwData{cursors}{cursors}
    \SetKwData{SignFlipVec}{SignFlipVec}
    \SetKwData{AnglesCol}{AnglesCol}
    \SetKwData{AnglesMtx}{AnglesMtx}
    \SetKwData{AnglesWork}{AnglesWork}
    \SetKwData{SortIndices}{SortIndices}
    \SetKwData{SubCursors}{SubCursors}
    \SetKwFunction{OrientEigenvectorsDecl}{\color{teal}{OrientEigenvectors}}
    \SetKwFunction{SortEigenvectors}{SortEigenvectors}
    \SetKwFunction{ReduceDimensionByOne}{ReduceDimensionByOne}
    \SetKwFunction{SolveRotationAnglesInSubDimDecl}{\color{teal}{SolveRotationAnglesInSubDim}}
    \SetKwFunction{SolveRotationAnglesInSubDim}{SolveRotationAnglesInSubDim}
    \SetKwFunction{ConstructSubspaceRotationMtx}{ConstructSubspaceRotationMtx}
    \SetKwFunction{MakeGivensRotation}{MakeGivensRotation}
    \SetKwProg{Def}{def}{:}{}
    %
    \Def{\OrientEigenvectorsDecl{\V, \E, \OrientToFirstOrthant}}{
      \Vsort, \Esort, \SortIndices $\leftarrow$ \SortEigenvectors{\V, \E}\;
      \Vwork $\leftarrow$ copy(\Vsort), $N$ $\leftarrow$ $\text{dim}(\V)$ \;
      \If{\OrientToFirstOrthant}{
        $\SignFlipVec[0]$ $\leftarrow$ $(\Vwork[0,0] => 0)$ ? $1$ : $-1$ \;
        $\Vwork[:, 0]$ $\leftarrow$ $\Vwork[:, 0] * \SignFlipVec[0]$ \;
      }
      \For{$i=0$ \KwTo $N-2$}{
        \Vwork, $\AnglesMtx[i, :]$ $\leftarrow$ \ReduceDimensionByOne{$N$, $i$, \Vwork} \;
      }
      $\SignFlipVec[-1]$ $\leftarrow$ $(\Vwork[-1,-1] => 0)$ ? $1$ : $-1$ \;
      \Vor $\leftarrow$ $\Vsort\cdot\diagoper(\SignFlipVec)$, \Eor $\leftarrow$ \Esort \;
      return \Vor, \Eor, \AnglesMtx, \SignFlipVec, \SortIndices   \;
    }
    \;
    %
    \Def{\SortEigenvectors{\V, \E}}{
      \SortIndices $\leftarrow$ argsort$(\text{abs}(\E))$ \;
      return $\Vsort[:, \SortIndices]$, $\Esort[\SortIndices]$, \SortIndices \;
     }
    \;
    %
    \Def{\ReduceDimensionByOne{$N$, $i$, \Vwork}}{
      \Vcol $\leftarrow$ $\Vwork[:, i]$ \;
      \AnglesCol $\leftarrow$ \SolveRotationAnglesInSubDim{$N$, $i$, \Vcol} \;
      \R $\leftarrow$ \ConstructSubspaceRotationMtx{$N$, $i$, \AnglesCol} \;
      \Vwork $\leftarrow$ $\R\transpose\Vwork$ \;
      return \Vwork, $\AnglesCol\transpose$ \;
    }
    \;
    %
    \Def{\SolveRotationAnglesInSubDimDecl{$N$, $i$, \Vcol}}{
      j $\leftarrow$ i + 1 \;
      \AnglesWork $\leftarrow$ zeros$(0: N)$ \;
      $\AnglesWork[j]$ $\leftarrow$
        arctan2$(\Vcol[j], \Vcol[j-1])$ \;
      \For{ $j = i + 2$ \KwTo $N$}{
        $\AnglesWork[j]$ $\leftarrow$
          arctan2$($ \;
          \hspace{3mm}$\Vcol[j] \bigl| \sin(\AnglesWork[j-1]) \bigr|, \bigl| \Vcol[j - 1] \bigr| $ \;
          $)$ \;
      }
      return $\AnglesWork$ \;
    }
    \;
    %
    \Def{\ConstructSubspaceRotationMtx{$N$, $i$, \AnglesCol}}{
      \R $\leftarrow$ $I(N)$ \;
      \For{ $j = N$ \KwTo $i+1$ \KwStep $-1$ }{
         \R $\leftarrow$ \;
         \hspace{3mm} $\MakeGivensRotation{$N, i, j, \AnglesCol[j]$} \cdot \R$  \;
      }
      return \R \;
    }
    \;
    %
    \Def{\MakeGivensRotation{$N, i, j, \theta$}}{
      \R $\leftarrow$ $I(N)$, $c$ $\leftarrow$ $\cos(\theta)$, $s$ $\leftarrow$ $\sin(\theta)$ \;
      \R$[i,i]$ $\leftarrow$ $c$, \R$[j,j]$ $\leftarrow$ $c$,
      \R$[i,j]$ $\leftarrow$ $s$, \R$[j,i]$ $\leftarrow$ $s$ \;
      return \R \;
    }
    \;
  }
  \caption{Pseudocode implementation of the \code{orient\_eigenvectors} function for the modified \code{arctan2} method in the absence of sparsity. Functions whose names are colored differ from the original algorithm. Matrix and vector indexing follows Python Numpy notation. Consult the source code in the \code{thucyd} package for the details around sparsity.
  }
\label{eigen::eq: algorithm pseudocode}
\end{algorithm}

There are two functions exposed on the interface of \code{thucyd.eigen}:
\begin{itemize}
  \item \code{orient\_eigenvectors} implements both the \code{arcsin} and modified \code{arctan2} methods. The pseudocode listing of Algorithm~\ref{eigen::eq: algorithm pseudocode} outlines the latter method. The function takes the keyword argument \code{method} that defaults to \code{arcsin} (for backward compatibility) but \code{arctan2} can be chosen and is now preferred.
  \item \code{generate\_oriented\_eigenvectors} is the same as detailed in the original article.
\end{itemize}
In certain application-specific scenarios, one may anticipate that the first eigenvector points into the first orthant, which is when all elements of the first eigenvector are positive. However, the sign returned from an \code{eig} or \code{svd} call may well be negative. Rather than accepting the major rotation that the modified \code{arctan2} method will impart, the \code{orient\_eigenvectors} function accepts an \code{OrientToFirstOrthant} flag that, when true, records a reflection at the top-level dimension, or \ebox{s_1 = -1}. The general case for the $S$ matrix when this flag is true is
\begin{equation}
  S_{\text{tan}} = \diagoper\left(\pm 1, 1, \ldots, 1, \pm 1 \right).
  \label{eigen::eq: S for arctan2 and orient to first orthant}
\end{equation}
With \ebox{s_1 = -1}, rotation $R_{1,2}\transpose$ is always minor.

%
%

\section{Empirical Example and Random Matrix Theory}
\label{eigen::sec: empirical example and random matrix theory}

\begin{figure*}[t]
  \centering
  \includegraphics[width=174mm]{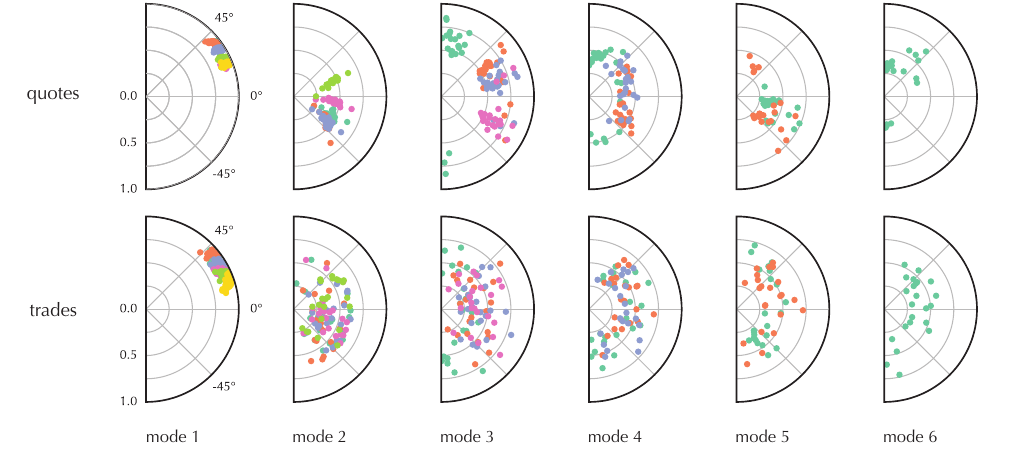}
  \caption{Hemispherical angles for quote and trade eigenvectors over 23 days and the six reducible modes. The angles were calculated using the \code{arcsin} method, and therefore they are all minor angles. The radii indicate the participation scores of the associated eigenvectors for each day: larger radii indicate better participation among the elements of an eigenvector.
  }
  \label{eigen::fig: eigenanalysis_quotes_and_trades_hemi_174}
\end{figure*}

\begin{figure*}[t]\sidecaption
  \centering
\includegraphics[width=127mm]{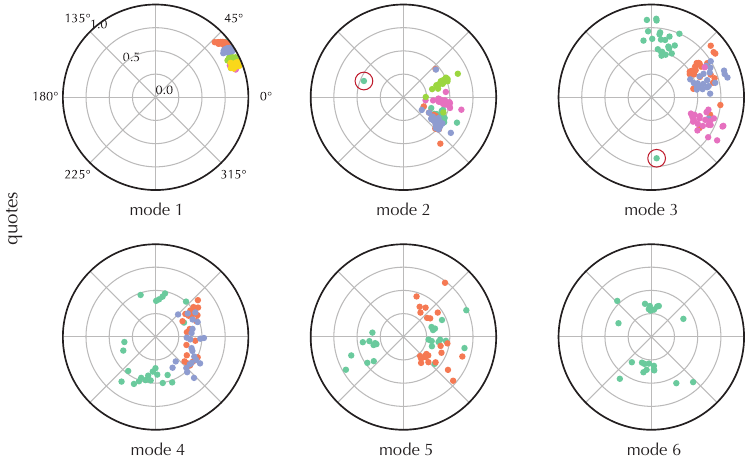}
\caption{Spherical angles of the quote eigenvectors over 23 days and the six reducible modes. The modified \code{arctan2} method was used to calculate the rotation angles. The first rotation angle of each mode spans the circle whereas subsequent angles are bound to the right hemisphere. Observe that the eigenvectors in mode 1 are highly directed whereas those in modes 4\textendash{}6 are well scattered. Modes 2 and 3 remain somewhat directed but with lower participation scores. The red circle highlights an apparent outlier that is discussed in the text.
}
\label{eigen::fig: eigenanalysis_quotes_sphr}
\end{figure*}

\begin{figure*}[t]\sidecaption
  \centering
  \includegraphics[width=127mm]{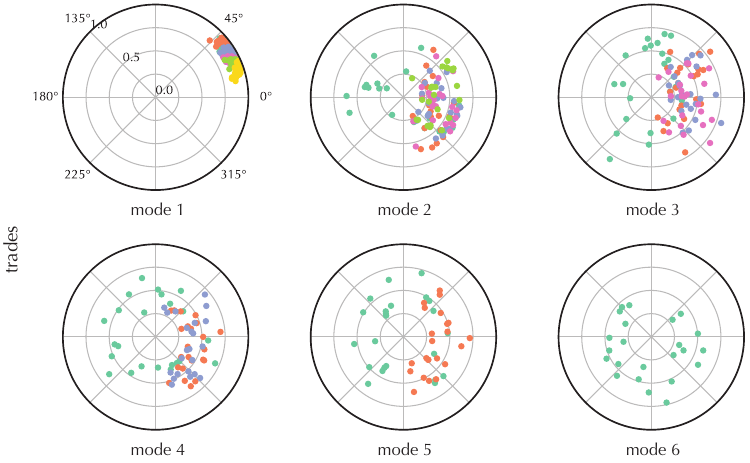}
  \caption{Spherical angles of the trade eigenvectors. Observe that the first mode is clearly directed into the first orthant whereas the other modes randomly scatter.
  }
  \label{eigen::fig: eigenanalysis_trades_sphr}
\end{figure*}

\begin{figure}[t]
  \centering
  \includegraphics[width=84mm]{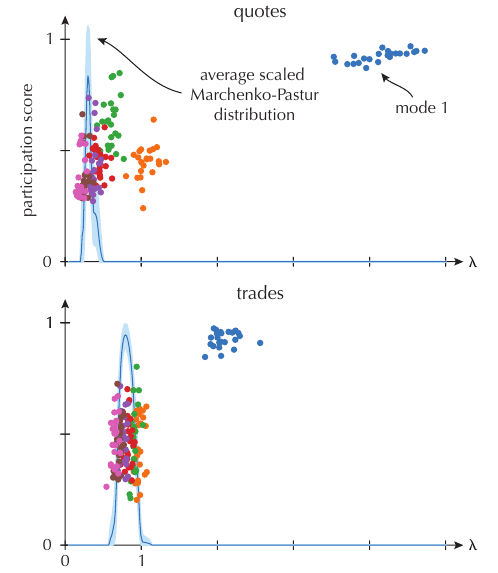}
  \caption{Overlay of empirical eigenvalue spectra with averaged \MarcenkoPastur distributions. (The colors here are used to distinguish the eigenmodes but do not correspond to the colors in the previous figures.) \emph{Top}, modes 1\textendash{}3 fall outside of the averaged MP distribution, indicating that their content has not been corrupted by sample noise.  \emph{Bottom}, only mode 1 falls well outside of the MP distribution while the rest lies firmly within the eigenvalue range that is indistinguishable from noise.
  }
  \label{eigen::fig: eigenanalysis_empirical_evals_vs_mp}
\end{figure}

The techniques of the previous sections are now applied to a dataset of quotes and trades from the financial markets in order to demonstrate the significance of unwrapping the $\pi$ interval of the \code{arcsin} method into the full $2\pi$ interval of the modified \code{arctan2} method. The consequent eigenvector pointing directions are then connected to the eigenvalue distribution of random matrix theory. The main results are 1) unwrapping the angular interval from $\pi$ to $2\pi$, as seen by comparing figure~\ref{eigen::fig: eigenanalysis_quotes_and_trades_hemi_174} with figures~\ref{eigen::fig: eigenanalysis_quotes_sphr}\textendash{}\ref{eigen::fig: eigenanalysis_trades_sphr}, allows for a distinction to be made between informative and uninformative eigenmodes; and 2) this distinction coincides well with what can be inferred from the eigenvalue spectra, as seen in figure~\ref{eigen::fig: eigenanalysis_empirical_evals_vs_mp}. These results then motivate static and dynamic methods to temporally stabilize eigenvectors, methods which are discussed in the following section.

The centered, standardized data panels~$P$ that are constructed below are analyzed via singular-value decomposition to produce
\begin{equation}
  P = U \Lambda^{1/2}\, V\transpose,
  \label{eigen::eq: svd def}
\end{equation}
where $U$ represents the projection of the data onto eigenvectors $V$ and $\Lambda$ is a diagonal matrix of eigenvalues. The eigendecomposition of the corresponding correlation matrix is
\begin{equation}
  \Corr(P) = V \Lambda\, V\transpose.
  \label{eigen::eq: eig def}
\end{equation}
At issue is the estimation of $V$ and $\Lambda$ from the data in $P$. The domain of classical statistics is asymptotically approached when the number of records $T$ for a fixed number of features $N$ goes to infinity: in this case, the Central Limit Theorem applies and estimation variance devolves to Gaussian. But in the real world, $T$ is finite and we can only hope to be in the \ebox{T > N} regime, for otherwise $P$ is singular. Sample noise \emph{always} exists with finite $T$, and so the question becomes one of distinguishing between real information and sample noise in $P$. The field of Random Matrix Theory (RMT) addresses this question \cite{Mehta:2004,Potters-Bouchaud:2021,Couillet:2022}.

The majority of the RMT literature until the early 2000s dealt with the probability density of the eigenvalues of a Gaussian random matrix, and within the universe of all random-matrix types, the real positive-definite matrix is relevant here. A principal result is that, in the limits \ebox{T\rightarrow \infty} and \ebox{N / T \rightarrow q > 0}, the eigenvalues $\lambda$ of such a random matrix follow the \MarcenkoPastur (MP) distribution
\begin{equation}
  \rho(\lambda) = \frac{ \sqrt{ (\lambda_+ - \lambda) (\lambda - \lambda_-) } }{ 2\pi q \lambda },
  \quad  \lambda_{\pm} = (1 \pm \sqrt{q})^2 ,
  \label{eigen::eq: MP pdf def}
\end{equation}
where \ebox{\lambda \in [\lambda_-, \lambda_+]}. Had $N$ been fixed, then $q\rightarrow 0$ and the density reduces to \ebox{\rho(\lambda) = \delta(\lambda - 1)}, which is to say, all eigenvalues are unity and so the eigenbasis is isotropic regardless of its orientation. When empirical eigenvalues fall outside of the MP distribution then they are not entirely corrupted by sample noise, but for those that fall within the distribution, the eigenvalues are indistinguishable from such noise.

As for eigenvectors $V$, it is easily established that the probability of drawing random matrix \ebox{E = P\transpose P / T} is invariant under the similarity transform \ebox{V' = O V O\transpose}, where $O$ is a pure rotation matrix. The consequence is that eigenvectors $v_k$ of such random matrices point uniformly across the surface of a hypersphere; this is known as the Haar measure. More recent advances in random matrix theory include the interaction between informative eigenvectors, known as rank-1 perturbations, to the eigenvectors that lie within the sample-noise bulk (see, for instance, \cite{Ledoit-Peche:2011,Bun-Bouchaud-Potters:2018}). However, research along this vein seldom speaks of adjusting the pointing direction(s) of informative eigenvectors (but see Bartz \etalwd \cite{Bartz:2011,Bartz:2016}); the main theme is how to better shrink empirical eigenvalues given the additional information available from the overlap between informative and randomly oriented eigenvectors. And so, there is an opportunity to use the techniques from section~\ref{eigen::sec: handedness in high dimensions} to distinguish between informative eigenvectors and those corrupted by noise, and to reorient the eigenvectors based on this distinction.

To begin the empirical example, Chicago Mercantile Exchange (CME) futures quote and trade data was obtained from the CME for seven foreign exchange contracts for May 2023. The FX pairs were EURUSD, USDJPY, GPBUSD, USDCHF, AUDUSD, NZDUSD, and CAD\-USD. The data covers 23~business days, some of which fall on regional holidays. The June 2023 contract was the most liquid during May and was therefore the focus. The daily datasets were sliced to retain records that fell between 00:00 and 16:00 Chicago time, which is generally the most liquid period. The quotes were for the best bid and offer and the trades were signed. A mid-price was generated from the quote data (ignoring changes in liquidity for the same price), and the mid-prices were processed to reduce the microstructure noise. Since the quote and trade data arrive asynchronously from contract to contract and from quote to trade, a common panel of records was created using time arbitration, one panel for each day. The panels had 14~columns, seven each for quotes and trades. The columns were then filtered with causal linear filters configured to a common timescale. The (log) quotes were filtered to measure price change while the (signed) trades were filtered to measure directionality. The panel was then downsampled in order to remove the autocorrelation imparted by the filters, leaving between 500\textendash{}1000 records for non-holiday days. Empirical study showed that each column was Student-$t$ distributed, so the dispersion and shape parameters were estimated, after which the columns were each mapped through the copula to a zero mean, unit variance Gaussian distribution, see \cite{Meucci:2007}. As a result, the final panel was multivariate Gaussian. This is the starting point for the eigenanalysis. Although the processing details have been only briefly summarized, they will be expanded upon in another publication.

While the panels were constructed with quotes and trades together in order to account for their relative update frequency, they can now be separated into two panels of 7~columns to inspect their inherent structure. Eigenanalysis was performed via SVD as in~(\ref{eigen::eq: svd def}) on the quote and trade panels, a separate analysis for each day. The rest of this section is devoted to the analysis of the eigenvalues~$\Lambda$ and eigenvectors~$V$.%
\footnote{%
See section~\ref{eigen::sec: declarations} for a link the $\Lambda$ and $V$ data used in this work.
} %
Next, the \code{orient\_eigenvectors} function was called with $V$ and $\Lambda$ using both \code{arcsin} and \code{arctan2} methods. Recall from equation~(21) of the original paper that the embedded angles can be organized into the upper triangle of a square matrix. For the panels in this example, the $\makevectorofsym{\theta}$ matrix looks like
\begin{equation}
  \begin{tikzpicture}[>=stealth,thick,baseline,
    every right delimiter/.append style={name=rd},
    ]
    {\small
      \matrix [matrix of math nodes,
          left delimiter=(,
          right delimiter=),
          ](A){
      0   & |[text=C1]|\theta_{1,2}
          & |[text=C2]|\theta_{1,3}
          & |[text=C3]|\theta_{1,4}
          & |[text=C4]|\theta_{1,5}
          & |[text=C5]|\theta_{1,6}
          & |[text=C6]|\theta_{1,7}  \\
          & 0
          & |[text=C1]|\theta_{2,3}
          & |[text=C2]|\theta_{2,4}
          & |[text=C3]|\theta_{2,5}
          & |[text=C4]|\theta_{2,6}
          & |[text=C5]|\theta_{2,7}  \\
          &
          & 0
          & |[text=C1]|\theta_{3,4}
          & |[text=C2]|\theta_{3,5}
          & |[text=C3]|\theta_{3,6}
          & |[text=C4]|\theta_{3,7}  \\
          &
          &
          & 0
          & |[text=C1]|\theta_{4,5}
          & |[text=C2]|\theta_{4,6}
          & |[text=C3]|\theta_{4,7}  \\
          &
          &
          &
          & 0
          & |[text=C1]|\theta_{5,6}
          & |[text=C2]|\theta_{5,7}  \\
          &
          &
          &
          &
          & 0
          & |[text=C1]|\theta_{6,7}  \\
          &
          &
          &
          &
          &
          & 0  \\
      };

      \draw[<-] (rd.east|-A-1-7) -- ++(0:5mm) node[right]{mode 1};
      \draw[<-] (rd.east|-A-6-7) -- ++(0:5mm) node[right]{mode 6};
    }
  \end{tikzpicture}
  \label{eigen::eq: theta-mtx in color}
\end{equation}
The diagonal color stripes highlight the common rotation order for the modes. For instance, the long green stripe indicates major angles, all of the other angles being minor. Now, there is a distinction between an eigenmode and a subspace that this matrix layout highlights: the modes of the system are represented by eigenvectors, which in turn populate the columns of $V$. Yet, the vectors are not independent from one another because they form an orthogonal basis. The eigenvector with the largest eigenvalue requires $N-1$ rotations to be oriented within $\realsspace^N$, after which the remaining vectors can only rotate in an orthogonal subspace. The next eigenvector then has only $N-2$ degrees of freedom for its orientation, and so on, until all reducible spaces are spanned. The eigenvector with the smallest eigenvalue lies in an irreducible space: it cannot be rotated because there are no axes that remain. As a consequence, the polar plots have fewer independent rotations (samples with different colors) as the angles for modes 1 through 6 are plotted.

%

To complement the angular information that appears in the plots, the radial coordinate plots the \emph{participation score} (PS). The participation score is a simple variation on the \emph{inverse participation ratio} (IPR) \cite{Monasson:2015}, defined as
\begin{equation}
  \text{IPR} = \sum_{i=1}^N v_i^4, \quad 1/N \leq \text{IPR} \leq 1,
  \label{eigen::eq: IPR def}
\end{equation}
where $v_i$ is the $i$th entry of an eigenvector in $V$, and \ebox{\makevector{v}\transpose\makevector{v} = 1} because $V$ is orthonormal. The lower bound of the IPR is attained when all entries $v_i$ participate equally, \ebox{v_i = 1/\sqrt{N}}, and the upper bound is attained when only one entry is nonzero, leaving \ebox{v_k = 1} for some $k$. In order to score higher participation above lower, the participation score goes as
\begin{equation}
  \text{PS} = 1 / (N \times \text{IPR}), \quad 1/N \leq \text{PS} \leq 1.
  \label{eigen::eq: PS def}
\end{equation}
When all elements participate equally, \ebox{\text{PS} = 1}. Now that the radial coordinate is defined we can move on to study the eigenvector systems.

Figure~\ref{eigen::fig: eigenanalysis_quotes_and_trades_hemi_174} plots the quote and trade angles on half-circles using on the \code{arcsin} method for the 23~days available. It is evident that mode~1 for quotes and trades is highly directed into the first orthant. Quote modes 2 and 3 appear somewhat directed, although with lower participation scores, and with a question mark on mode 3 because of the appearance of samples toward the south pole of the plot: have these been wrapped around from the north? Trade modes 2\textendash{}6 have their angles fully scattered across $\pi$,  and the participation scores are small and rather unstable. Quote mode 4 has randomly scattered angles and persistently low scores; it is hard to tell what to say about modes 5 and 6.

Figures~\ref{eigen::fig: eigenanalysis_quotes_sphr} and~\ref{eigen::fig: eigenanalysis_trades_sphr} replot the quote and trade eigensystems after the angles were calculated using the modified \code{arctan2} method. The green points have been expanded into major angles. It is now evident that quote mode~4 truly points in random directions, and that is probably the case for modes 5 and 6. The full-circle plots for modes 2 and 3 are also more revealing: the modes are indeed directed, albeit less so. For mode 3, the points that had appeared in the south now all point to the north except for an outlier, and this highlights the benefit of moving from the \code{arcsin} method to the \code{arctan2} method. The outlier, which is encircled in the graphs for modes 2 and 3, can only be discerned after switching to the \code{arctan2} method. This outlier is not an artifact of the orientation algorithm but is an actual outlier in the data: the monthly update of the US Consumer Price Index (CPI) was released on May 10, 2023. Market participants were following the CPI closely in order to assess whether the US Federal Reserve would curtail their rate hikes. There was a price jump at the time of the announcement for all of the FX pairs considered, but EUR\-USD responded anomalously. (Removal of this pair from the panel construction restored the outlining point to the north-facing cluster.) In practice, the quantitative analyst would likely split the data into times outside of scheduled economic events from those times around the events, and then analyze the split-out data separately. In this way, ``normal'' market behavior is separated from times when the markets often jump to new levels, but for the purpose here it is more interesting to show that the \code{arctan2} can reveal an outlier.

The trade eigensystems shown in figure~\ref{eigen::fig: eigenanalysis_trades_sphr} show how the major angle scatters uniformly over the full circle for modes 2\textendash{}6 while the minor angles continue to scatter in the right hemisphere. Only mode 1 is highly directed and exhibits high participation scores.

These observations can now be compared to the traditional RMT eigenvalue analysis. Figure~\ref{eigen::fig: eigenanalysis_empirical_evals_vs_mp} plots empirical eigenvalues and participation scores as points on a Cartesian graph, and plots an average of the associated rescaled \MarcenkoPastur distributions. As discussed above, almost any analysis of empirical data panels is corrupted by noise induced by the finite length of the panels. For instance, eigenmodes whose eigenvalues fall within the MP distribution, while produced from real data, are indistinguishable from noise. The \emph{top} and \emph{bottom} panes show that quote modes 1\textendash{}3 and trade mode 1 fall outside of the (average) MP distribution; the other modes fall within the distribution. This is a very good match to the conclusions made above.

However, this eigenvalue analysis is not absolute, and it is here where the analysis of pointing directions can be leveraged: a rescaling of the MP distribution was performed, and while there is no question that rescaling is necessary, the results can be somewhat self-fulfilling. This is what is going on: the MP distribution applies to purely random matrices, but here we have a mixture of informative and random modes. For trades, there is only one informative mode, so really there are $N-1$ modes that are indistinguishable from noise. The $q$ value needs to be adjusted to read \ebox{q' = (N - 1) / T}; this change is simple. In addition, eigenvalue $\lambda_1$ has contributed to the scaling of the multivariate Gaussian distribution, so this needs correction \cite{Potters-Bouchaud:2021,Meucci:Labs:RMTShrinkage}. With
\begin{equation}
  \bar{\lambda} = \text{Mean} \left( \lambda_2, \lambda_3, \ldots, \lambda_N \right),
  \label{eigen::eq: lambda-bar MP rescaling}
\end{equation}
the rescaled MP distribution reads
\begin{equation}
  \rho_{\bar{\lambda}}(\lambda) = 1 / \bar{\lambda}\; \rho\left( \lambda / \bar{\lambda}  \right).
  \label{eigen::eq: MP pdf rescaled}
\end{equation}
It is the average (over the 23 days) of these rescaled distributions that are plotted in the trades pane of the figure. For quotes, three modes were deemed to be informative, leaving four modes to generate noise (the irreducible mode has to be included). Every time an eigenmode is excluded from the noise, the remaining distribution shifts to the left, leading to a potentially subjective line between modes to include or exclude. It is here where the analysis of pointing directions introduces new information to cross-validate the traditional RMT eigenvalue approach.

%
%

\section{Eigenvector Stabilization}
\label{eigen::sec: Eigenvector Stabilization}

The empirical eigenvector pointing directions studied in the last section exhibit two different ways in which they evolve in time. Informative modes are at least somewhat directed, and so their average direction is meaningful. Noninformative modes, in contrast, scatter randomly on the hypersphere, and so their average, which is not meaningful, can be replaced by a fixed direction, thereby exchanging variance for bias. Let us call these two methods of stabilization \emph{dynamic} and \emph{static}. These methods could be used in future work to help clean empirical correlation matrices. For now, figures~\ref{eigen::fig: eigenanalysis_correlation_grids} and~\ref{eigen::fig: eigenanalysis_correlation_dispersions} show examples of how the stabilization of the eigenvectors in turn stabilizes the evolution of the correlation matrix. But first, these methods are quantified.

\begin{figure}[t]
  \centering
  \includegraphics[width=84mm]{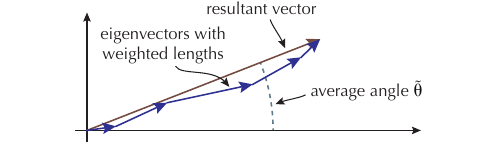}
  \caption{The average pointing direction in $\realsspace^2$ over five wobbling vectors is found by stacking the vectors end-to-end and then measuring the resultant vector. Applying weights to the vectors is equivalent to scaling the vector lengths, as illustrated here.
  }
  \label{eigen::fig: eigenanalysis_average_angle}
\end{figure}
%

\begin{figure*}[t]
  \centering
  \includegraphics[width=174mm]{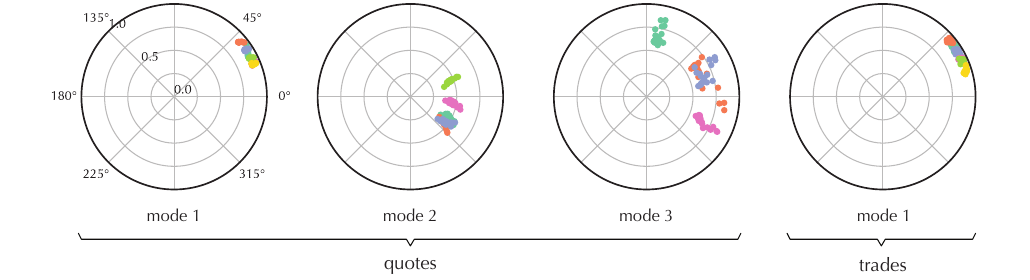}
  \caption{Angular distributions after five-point filtering of the directed quote and trade eigenvectors. To the \emph{left} are the angles of quote modes 1\textendash{}3, and to the \emph{right} lie the angular distributions of the trade mode~1. The filter weights were \ebox{h[n] = [1, 2, 3, 2, 1] / 9}, which is the result of convolving a three-point box filter with itself. The average delay for this $h[n]$ is two samples, or two days here. Observe the stabilization of the eigenvector wobble and the concentration of the participation scores, especially for quote and trade modes 1.
  }
  \label{eigen::fig: eigenanalysis_average_pointing_directions}
\end{figure*}

With respect to dynamic stabilization, it is evident from figures~\ref{eigen::fig: eigenanalysis_quotes_sphr} and~\ref{eigen::fig: eigenanalysis_trades_sphr} that the pointing directions of the informative modes wobble in time, whereas, in contrast, the eigenvectors of the random modes scatter uniformly. A causal local average of the latest pointing directions can be estimated by applying a simple filter, but how shall this be performed? Recall from section~9.2 of the original article that the angular inclination of a collection of vectors is correctly computed by stacking the vectors together end-to-end and then calculating the angle of the resultant vector, such as pictured in figure~\ref{eigen::fig: eigenanalysis_average_angle}. This stands in contrast to simply averaging the angles directly. In analogy to equation~(31) in the original article, we can write more generally
\begin{equation}
  M_h[n] = \left( h \conv \orientedeigvec  \right) [n]
  \label{eigen::eq: filtered eigenvector mtx}
\end{equation}
and then renormalize via
\begin{equation}
  M[n] = M_h[n] \; \ltwonorm{M_h[n]}\inverse
  \label{eigen::eq: Mh renormalized}
\end{equation}
where
\begin{equation*}
  \ltwonorm{M_h}
  \equiv
  \diagoper \Bigl( \ltwonorm{M_{h,1}} , \ltwonorm{M_{h,2}}, \ldots, \ltwonorm{M_{h,N}} \Bigr).
\end{equation*}
What is happening here is that a time series of orientated eigenbases $\orientedeigvec[n]$ is convolved with impulse response $h[n]$ to yield $M_h[n]$, and to make sense in this context, the weights of $h[n]$ are all positive. For any~$n$, $M_h$ is no longer orthonormal, for this property is not guaranteed under basis addition. Subsequent normalization into $M[n]$ is done in~(\ref{eigen::eq: Mh renormalized}). To orthogonalize, the orientation algorithm can be used by computing the sequence
\begin{equation}
  M[n] \longrightarrow \makevectorofsym{\theta}_M[n] \longrightarrow \bar{\orientedeigvec}[n],
  \label{eigen::eq: M -> V-bar orthogonalization}
\end{equation}
where calls to \code{orient\_eigenvectors} and \code{generate\_} \code{oriented\_eigenvectors} are using in turn. From these intermediate steps we can ultimately write
\begin{equation}
  \left( h \conv \orientedeigvec \right) \! [n] \longrightarrow \bar{\orientedeigvec}[n].
  \label{eigen::eq: (h conv V) -> V-bar}
\end{equation}
Let us denote by \ebox{\tilde{\makevectorofsym{\theta}}[n]} the matrix of embedded angles associated with \ebox{\bar{\orientedeigvec}[n]}, where the tilde on $\makevectorofsym{\theta}$, rather than a bar, emphasizes that the average here is not over angles but eigenvectors.

For consistency with the eigenvector treatment, the eigenvalues are also filtered; we have
\begin{equation}
  \bar{\Lambda}[n] = \left( h \conv \Lambda \right) \! [n].
  \label{eigen::eq: (h conv L) -> L-bar}
\end{equation}
The rank order of eigenvalues does not change upon convolution with positive $h[n]$, so there is no disruption in the map between eigenvalues to eigenvectors.

Figure~\ref{eigen::fig: eigenanalysis_average_pointing_directions} shows the result of filtering the quote and trade eigenvectors for the informative modes. A five-point causal filter was applied, so each point represents the local average of five contiguous days. The caption covers the filter details. A statistical analysis would require the downsampling of the filtered points, but here all filtered points were plotted because there are so few points in this dataset. As expected, the scatter of filtered pointing directions is reduced. In addition, the participation scores have tightened and slightly increased. Quote mode 2 continues to have a curiously low participation score, but the vectors remain directed. All in all, eigenvector wobble has been demonstrably reduced by applying a filter. Increasing the filter length would add further stability but at the expense of adding more delay into the system.

With respect to static stabilization, we have already seen that quote modes 3\textendash{}6 and trade modes 2\textendash{}6 scatter randomly, so taking their local average is meaningless. An evident way to statically reorient the eigenvectors is to fix the angles of the random modes, and in the absence of a better alternative, the angles can simply be set to zero. In doing so, an angular bias is put in place in order to reduce out-of-sample error. For the quotes, angle matrix $\makevectorofsym{\theta}$ that appears in~(\ref{eigen::eq: theta-mtx in color}) is rewritten as
\begin{equation}
  \makevectorofsym{\theta}_{\text{modal}} =
  \left(
    \begin{array}{ccccccc}
      0
        & \theta_{1, 2}
        & \theta_{1, 3}
        & \theta_{1, 4}
        & \theta_{1, 5}
        & \theta_{1, 6}
        & \theta_{1, 7}  \\
        & 0
        & \theta_{2, 3}
        & \theta_{2, 4}
        & \theta_{2, 5}
        & \theta_{2, 6}
        & \theta_{2, 7}  \\
        &
        & 0
        & \theta_{3, 4}
        & \theta_{3, 5}
        & \theta_{3, 6}
        & \theta_{3, 7}  \\
        &  &  & 0  & 0  & 0  & 0  \\
        &  &  &    & 0  & 0  & 0  \\
        &  &  &    &    & 0  & 0  \\
        &  &  &    &    &    & 0  \\
    \end{array}
  \right),
  \label{eigen::eq: theta-mtx statically stabilized}
\end{equation}
where $\makevectorofsym{\theta}_{\text{modal}}$ only retains nonrandom modes. Static stabilization can be applied to either daily oriented matrices $\orientedeigvec[n]$ or their averages $\bar{\orientedeigvec}[n]$. A statically stabilized matrix $V$ is then reconstructed by calling \code{generate\_} \code{oriented\_eigenvectors} with $\makevectorofsym{\theta}_{\text{modal}}$. Static stabilization has the flavor of Principal Component Analysis (PCA) but it is not. With PCA, informative eigenmodes are retained while the others are discarded, but with static stabilization, all modes remain present in $\orientedeigvec$. Geometrically, what is going on is that $\orientedeigvec_{\text{modal}}$ is ``closer'' to identity matrix $\Identity$ than $\orientedeigvec$ itself. For the six reducible modes in this example, expanding~(\ref{eigen::eq: V-orient = V S = R}) leads to
\begin{equation}
  \orientedeigvec_{\text{modal}} = R_1 R_2 R_3 R_4 R_5 R_6 \Identity = R_1 R_2 R_3 \Identity
  \label{eigen::eq: V-modal via rotation}
\end{equation}
since we have $R_{4,5,6} = \Identity$.

We are now in a position to study the correlation structure as the stabilization steps are applied. A correlation matrix is calculated according to
\begin{equation}
  \overline{\Corr}(P) =
    \diagoper(\makevectorofsym{s}\inverse)
    \left(
      \bar{\orientedeigvec} \bar{\Lambda} \bar{\orientedeigvec}\transpose
    \right),
    \diagoper(\makevectorofsym{s}\inverse),
  \label{eigen::eq: corr from V and L}
\end{equation}
where \ebox{\makevector{s}^2 = \diagoper\left( \bar{\orientedeigvec} \bar{\Lambda} \bar{\orientedeigvec}\transpose \right)}. Note that while \ebox{\orientedeigvec\Lambda\orientedeigvec\transpose} produces a proper correlation matrix, being positive definite with ones along the diagonal, \ebox{\bar{\orientedeigvec} \bar{\Lambda} \bar{\orientedeigvec}\transpose} does not. The latter product is positive definite, but the diagonal entries have small deviations from unity. Such imperfections arise from the processing used to produce $\bar{\orientedeigvec}$ and $\bar{\Lambda}$, and are removed by the two end terms in~(\ref{eigen::eq: corr from V and L}).

\begin{figure}[t]
  \centering
  \includegraphics[width=84mm]{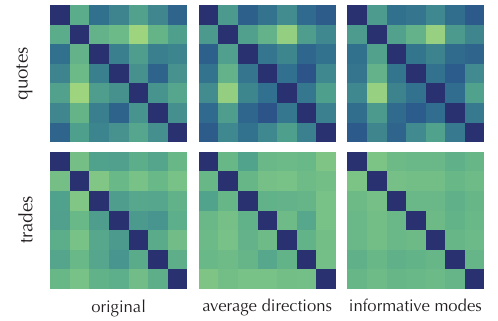}
  \caption{Quote and trade correlation matrices, raw through fully stabilized, for a single day. The original correlation matrices (\emph{left}) contain a lot of structure. The dynamically stabilized correlations (\emph{center}) have less structure. Lastly, the dynamically and statically stabilized correlations (\emph{right}) have the least structure and are the most persistent over time.
  }
  \label{eigen::fig: eigenanalysis_correlation_grids}
\end{figure}

Figure~\ref{eigen::fig: eigenanalysis_correlation_grids} shows a progression of matrix stabilization from the original to dynamically stabilized to both dynamically and statically stabilized. The correlation matrices in the figure are for one day, a single day for the original, and a single five-day average for the dynamically filtered. What is apparent as the eigenvectors are stabilized is that structure in the correlation matrix is removed; this is especially true for the trade correlations. Although the effect of static stabilization is somewhat apparent for this progression of matrices, the benefit of filtering the eigenvector directions cannot be shown from a single snapshot.

\begin{figure*}[t]
  \centering
  \includegraphics[width=174mm]{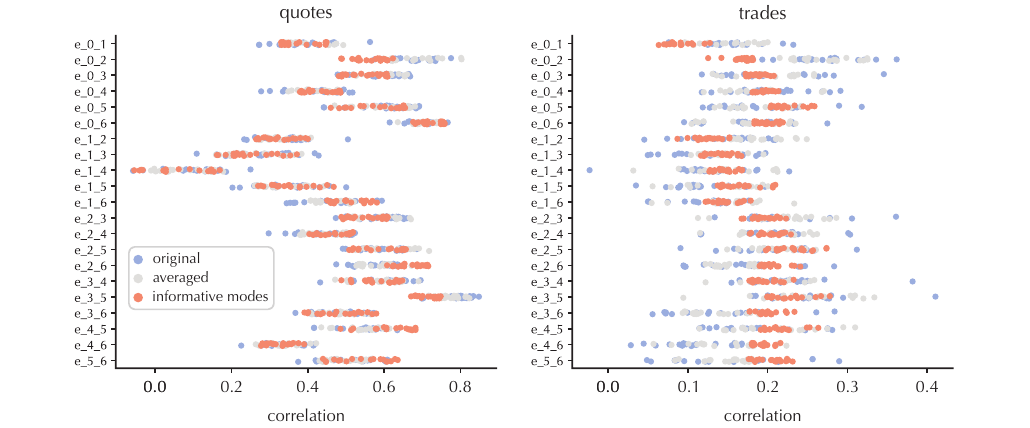}
  \caption{Pairwise correlation scatter for each element of the correlation matrices in figure~\ref{eigen::fig: eigenanalysis_correlation_grids}. Observe the reduction in correlation dispersion as the eigensystem is stabilized, first by filtering the eigenvector pointing directions (gray circles) followed by retaining the angular structure of only the informative modes (orange circles).
  }
  \label{eigen::fig: eigenanalysis_correlation_dispersions}
\end{figure*}

To better investigate the effect of the stabilization, the correlation scatter of each entry in the correlation matrices can be plotted over the population, which in this case is 23~days less the 5~day filter window, or 18 available samples. This analysis appears in figure~\ref{eigen::fig: eigenanalysis_correlation_dispersions}. Overall, correlations (gray circles) that result from filtering the eigenvectors to create local averages fall with\-in the range of the original correlations (blue circles). This is the consequence of wobble reduction. However, the randomly oriented modes still have an effect. The second stabilizing stage, where only the angles of the informative modes are retained, further decreases the dispersion of the correlations (orange circles).

Static stabilization brings us close to the topic of shrinkage of a correlation matrix subject to corruption from finite-sample noise. Recall that \LedoitWolf shrinkage, being the matrix generalization of \JamesStein shrinkage, blends an empirical matrix with an identity matrix of the same size,
\begin{equation}
  \Sigma_{\text{shr}} = \alpha \makeestimate{\Sigma} + (1 - \alpha) \Identity,
  \label{eigen::eq: LW shrinkage}
\end{equation}
where \ebox{\alpha \in [0, 1]} is the shrinkage factor \cite{Ledoit-Wolf:2004}. However, this estimator applies to empirical correlation matrices that have no informative modes (often called ``spikes'' in the literature; see, for instance, \cite{Benaych:2011}). Informative modes cannot simply be dropped from a matrix of eigenvectors (in the manner that PCA drops uninformative modes) because the remaining vectors still include all $N$ entries, but we \emph{can} rotate them away in order to isolate the remaining uninformative modes. With the quote panel in mind, we figuratively have
\begin{equation}
  \makeestimate{\Sigma}_{R} =
  R_3\transpose R_2\transpose R_1\transpose \orientedeigvec =
  \left(
    \begin{array}{ccccccc}
      1  &    &    &             &             &             & \\
         & 1  &    &             &             &             & \\
         &    & 1  &             &             &             & \\
         &    &    & \centerdot  & \centerdot  & \centerdot  & \centerdot  \\
         &    &    & \centerdot  & \centerdot  & \centerdot  & \centerdot  \\
         &    &    & \centerdot  & \centerdot  & \centerdot  & \centerdot  \\
         &    &    & \centerdot  & \centerdot  & \centerdot  & \centerdot
    \end{array}
  \right),
  \;
\end{equation}
which is admissible into~(\ref{eigen::eq: LW shrinkage}), after which the shrunk matrix is rotated back. Only the noise-related eigenvalues should be included when calculating $\alpha$. Tying \LedoitWolf shrinkage back to static stabilization, the effective shrinkage weight for the latter is \ebox{\alpha = 0}.

To summarize this section, directional analysis has been applied twice to the eigenanalysis of quote and trade data as it evolves in time. In the first instance, the pointing directions of the eigenvectors have been filtered to reduce wobble. In the second instance, the orientation of randomly oriented modes was stabilized by setting their embedded angles to zero in order to trade off out-of-sample error with bias. The effects of this type of matrix ``cleaning'' were demonstrated by taking several types of views of the data, from angular to correlation to pairwise correlation dispersion. When a correlation matrix contains a mixture of informative and noise-corrupted modes, one can rotate away the informative modes, apply proper shrinkage, and then rotate the system back.

%
%

\section{Conclusions}
\label{eigen::sec: conclusions}

Two methods now exist to solve the system of transcendental equations necessary to rotate an eigenbasis into the identity matrix: In the original article, which introduced the \code{arcsin} method, rotations and reflections were combined to scan for, and rectify, reflections that could be embedded in an eigenvector matrix $V$ returned from software calls to \code{svd} or \code{eig}. This article introduced the modified \code{arctan2} method that can defer the rectification of reflections embedded in $V$ until the last subspace. In doing so, \ebox{N-1} angles within \ebox{V \in \realsspace^{N\times N}} can span a $2\pi$ interval rather than the single $\pi$ interval of the \code{arcsin} method. This in turn allows for the disambiguation of angular wrap-around on the $\pi$ interval for these rotations, and the empirical example showed the importance of observing the full sweep of eigenvector pointing directions in an evolving system. The \code{arcsin} and modified \code{arctan2} methods serve different purposes: the former is better for uses such as regression where the instability of eigenvector signs along an evolving datastream corrupts downstream interpretability, and the latter is better for determining the degree of eigenvector directedness over time.


\begin{acknowledgements}

I am grateful for my conversations with Dr. A. Meucci about the distribution mapping through the copula; this led to a cleaner application of random matrix theory.

\end{acknowledgements}


%
%

\section{Declarations}
\label{eigen::sec: declarations}

\paragraph{Conflict of interest:} I declare that I have no conflict of interest.

\paragraph{Code availability:} The Python package \code{thucyd} is available on \textsc{PyPi} and \textsc{Conda-Forge} at
\begin{itemize}
  \item \href{https://pypi.org/project/thucyd/}{\code{pypi.org/project/thucyd}}
  \item \href{https://github.com/conda-forge/thucyd-feedstock}{\code{github.com/conda-forge/thucyd-feedstock}} .
\end{itemize}
This paper covers \code{thucyd} version 0.2.5+. The code is covered by the \textsc{Apache 2.0} license and is freely available. In addition, the source code for the reference implementation is freely available at
\begin{itemize}
  \item \href{https://gitlab.com/thucyd-dev/thucyd}{\code{gitlab.com/thucyd-dev/thucyd}} .
\end{itemize}

\paragraph{Data availability:} The CME data used to generate the panel data cannot be republished. However, the eigensystems, which are the focus of this paper, are available at
\begin{itemize}
  \item \href{https://gitlab.com/thucyd-dev/thucyd-eigen-working-examples}{\code{thucyd-eigen-working-examples}} .
\end{itemize}
Follow the project's \code{Readme.md} file for details.


\bibliographystyle{spmpsci}
\bibliography{arxiv_eigen_biblio}

\typeout{get arXiv to do 4 passes: Label(s) may have changed. Rerun}

\end{document}